\documentclass{ws-ijgmmp}
\usepackage{dsfont}
\usepackage{color}
\def\bb#1\eb{\textcolor{blue}
{#1}} %
\usepackage[normalem]{ulem}
\begin{document}

\markboth{Miguel A. Javaloyes And Miguel S\'anchez}
{Finsler metrics and relativistic spacetimes}

%
\catchline{}{}{}{}{}
%

\newcommand{\R}{\mathds R}
\newcommand{\N}{\mathds N}
\newcommand{\Z}{\mathds Z}

\title{Finsler metrics and relativistic spacetimes}

\author{Miguel A. Javaloyes}

\address{Departamento de Matem\'aticas, Universidad de Murcia,\\ Facultad de Matem\'aticas, Campus de Espinardo s/n\\
Murcia, E-30100, Spain\,\\
\email{majava@um.es} }

\author{Miguel S\'anchez}

\address{Departamento de Geometr\'{\i}a y Topolog\'{\i}a, Universidad de Granada.\\ Facultad de Ciencias, Campus de Fuentenueva s/n. \\ E-18071 Granada (Spain)
\\
\email{sanchezm@ugr.es} }

\maketitle

\begin{history}
\received{(Day Month Year)}
\revised{(Day Month Year)}
\end{history}

{\small {\em Submitted to the Special Issue for the XXII IFWGP
Evora, with Associated Editors R. Albuquerque, M. de Leon and M.
C. Munoz Lecanda.}}

\bigskip

\begin{abstract}
Recent  links between Finsler Geometry and the geometry of
spacetimes are briefly revisited,  and prospective ideas and
results are explained. Special attention is paid to geometric
problems with a direct motivation in Relativity and other parts of
Physics.
\end{abstract}


\section{Introduction}

There has been a recent interest in the links between Finsler
Geometry and the geometry of relativistic spacetimes. A well-known
reason comes from the viewpoint of  applicability: Finsler
metrics are much more general than Riemannian ones and,
accordingly, if one replaces the Lorentz metric of a spacetime by
a Finslerian counterpart, the possibility to model physical
effects is richer. A more practical reason comes from a purely
mathematical correspondence: the geometry of a concrete class of
Finsler manifolds (Randers spaces) is closely related to the
conformal structure of a  class  of spacetimes (standard
stationary ones). So, results in one of these two fields can be
translated into results on the other one ---and, sometimes, this
can be generalized to more general  Finsler and Lorentz spaces.
The purpose of this paper is to  make a brief revision of this
topic, emphasizing the  applications to Physics.

After some preliminary definitions and motivations on Finsler
metrics (Section \ref{s2}), we focus on the problem of describing
spacetimes by using Finsler elements (Section \ref{s3}).  We
introduce and discuss a notion of Lorentz-Finsler metric which
leads to {\em (conic) Finsler spacetimes}. Our definition is quite
general and, so, several previous notions in the literature fit in
it. Such a Finsler spacetime provides a {\em cone structure}. This
notion  has been studied systematically very recently \cite{FS},
and it  allows  us to introduce classical Causality in the
framework of Finsler spacetimes. In particular, some of the key
results on this topic can be recovered, Theorem \ref{tFS}. It is
worth pointing out: (a) the cone structure of a Finsler spacetime
generalizes the one provided by the chronological futures of a
classical spacetime, and defines implicitly a past cone structure
(but no further assumption on reversibility should be necessary if
one  liked to include in the definition of conic Finsler spacetime
the two cones at each point), and (b) in classical spacetimes, the
cone structure is equivalent to the conformal one, but in Finsler
spacetimes many non-conformally related Lorentz Finsler metrics
will have the same cone structure and, thus, the same Causality
---that is, the information of the Lorentz-Finsler metric not
contained in the cone structure is much richer  than in classical
spacetimes.

In Section \ref{s4}, we make a brief survey on the recent
 developments   about the geometric correspondence between standard
stationary spacetimes and Randers spaces. We focus in those
questions with more clear applications   to relativistic
spacetimes and, so, topics such as the causal structure of
standard stationary spacetimes or the gravitational lensing are
emphasized. Finally, the conclusions are summarized in the last section.

\section{Finsler metrics and their applications } \label{s2}

Let us introduce the very general notion of (conic) pseudo-Finsler
metric in a manifold $M$. Let  $TM$ denote the tangent bundle of
$M$, $\pi:TM\rightarrow M$, the natural projection and $A\subset
TM\setminus \bf 0$,  an  open  subset satisfying that it is conic
(that is, if $v\in A$ and $\lambda>0$, then $\lambda v\in A$), and
it projects on all $M$ (i.e., $\pi(A)=M$). We say that a  smooth
function $L:A\subset TM\setminus {\bf 0}\rightarrow \R$ is a {\it
(conic) pseudo-Finsler metric} if $L$ is positive homogeneous of
degree 2, namely, $L(\lambda v)=\lambda^2 L(v)$ for any $v\in A$
and $\lambda>0$, and the fundamental tensor $g$, defined as the
Hessian of $\frac 12 L$  at every $v\in A$, is nondegenerate. In
other words, given $v\in A$, the bilinear symmetric form in
$T_{\pi(v)}M$ defined as
\begin{equation}\label{fundten}
g_v(u,w)=\frac{1}{2}\left. \frac{\partial^2}{\partial s\partial t}L(v+tu+sw)\right|_{t=s=0}
\end{equation}
for any $u,w\in T_{\pi(v)}M$, is non-degenerate.

Let us remark that classical Finsler metrics are a particular case
of pseudo-Finsler metrics. More specifically, $L$ is a Finsler
metric if $A=TM\setminus \bf 0$ and  the fundamental tensor $g$ is
positive definite\footnote{In this case, positive
definiteness follows from non-degeneracy \cite[Prop.   2.16
(ii)]{JS}.} for every $v\in A$. In this case, $L$ is always
positive, since by homogeneity $L(v)=g_v(v,v)$ for every $v\in A$
and we can consider $F=\sqrt{L}$, which is positive homogeneous of
degree one. The function $F$ is usually called the Finsler metric
function,  and it is determined by its {\em indicatrix}, i.e. the
set of its unit vectors.

 Among the most classical examples of Finsler metrics,  we will use Randers metrics, which are given by
\begin{equation}\label{randers}
F(v)=\alpha(v)+\beta(v)
\end{equation}
for every $v\in TM$, where $\alpha(v)=\sqrt{h(v,v)}$ and $h$ and
$\beta$ are a Riemannian metric and a one-form in $M$
respectively. In fact, the fundamental tensor of $F$ is positive
definite in $v\in TM$ if and only if $\alpha(v)+\beta(v)>0$ (see
\cite[Corollary 4.17]{JS});   this restriction is satisfied
everywhere if the $h$-norm of $\beta$ is smaller than 1, so
defining a  Randers metric. These metrics appear naturally in
several contexts as in Zermelo navigation problem, which aims to
describe the trajectories that minimize the time in the presence
of a mild wind or current modelled by a vector field $W$. Then,
these trajectories are given by geodesics of a Randers metric
defined as
\begin{equation}\label{zermelo}
F(v)=\frac{-g(v,W)+\sqrt{g(v,W)^2+g(v,v)(1-g(W,W))}}{1-g(W,W)}
\end{equation}
where $v\in TM$, $g$ is a Riemannian metric and $g(W,W)<1$ in all
$M$ \cite{BRS04}. Observe that Zermelo metrics are always positive
for every $v\in TM\setminus \bf 0$ and then its fundamental tensor
is positive definite. As we will see later, Randers metrics also
appear naturally associated to stationary spacetimes describing
their causal properties (see Section \ref{s4}). Other remarkable
example is given by Matsumoto metric, which describes trajectories
minimizing time in the presence of a slope  ---recall that
going up is slower than going down. This metric is defined as
\[F(v)=\frac{\alpha(v)^2}{\alpha(v)-\beta(v)}\]
for every $v\in TM$, and its fundamental tensor is positive
definite  in $v\in TM\setminus\bf 0$ if and only if
$(\alpha(v)-\beta(v))(\alpha(v)-2\beta(v))>0$  \cite[Corollary
4.15]{JS}.  Therefore, Matsumoto metric is properly {\em conic} if
this inequality is not satisfied by some $v\neq 0$. Randers and
Matsumoto metrics are particular examples of the class of
$(\alpha,\beta)$-metrics, which are defined as $F(v)=\alpha
\phi(\beta(v)/\alpha(v))$, being $\phi$ an arbitrary non-negative
real function  (see \cite{JS} and references therein).

\section{Finsler-Lorentz metrics and spacetimes}\label{s3}
In classical General Relativity, a {\em spacetime} is a
(connected, Hausdorff) $n$-manifold $M$ endowed with a Lorentzian
metric $g$ and a time orientation, i.e., a continuous choice of
one of the two timelike cones at each point, which will be
regarded as {\em future-directed}. There are some speculative
applications of the replacement of $g$ by a   (generalization of
a) Finsler metric $F$, as  modeling possible anisotropies of the
spacetime even at an infinitesimal level, or admitting speeds
higher than light. In any case, the possibility to model a general
action functional (homogeneous of degree two, but not necessarily
coming from a quadratic form) justifies the study of the
Lorentz-Finsler approach.

There are different possibilities in order to define a Finsler spacetime. Recall that, for a Riemannian metric $g_R$, the unit vectors constitute the indicatrix of a (standard) Finsler metric; in particular, each unit sphere $S_p\subset T_pM$  is convex (as the boundary of the unit ball) at each point. Nevertheless, if $g$ is a Lorentzian metric, one has to consider the subsets $S_p^+, S_p^-\subset  T_pM$ containing, resp. the spacelike and timelike unit vectors at $p$. Notice that
$S_p^-$ always contains two connected parts, each one {\em concave}. For $n=2$ these properties also hold for $S_p^+$, but for higher dimensions $S_p^+$   is connected and, at each point $p$, its second fundamental form (say, with respect to any auxiliary Euclidean product at $T_pM$) has Lorentzian signature.  From the physical viewpoint, the most important elements of the spacetime are the causal vectors of the metric, since they model the trajectories of massive and massless particles. Observe in particular that the length of a causal curve $\gamma:[0,1]\rightarrow M$ in the spacetime is computed as
\[\ell_g(\gamma)=\int_0^1 \sqrt{-g(\dot\gamma,\dot\gamma)}ds\]
and a fundamental property in the spacetime is that causal
geodesics locally maximize this length. If we want to define a
more general way of measuring the length of curves  (but
preserving that the length does not depend on the parametrization
of the curve and geodesics are local length-maximizers), then we
must consider a positive one-homogeneous function $F:A\subset
TM\setminus {\bf 0}\rightarrow (0,+\infty)$ with fundamental
tensor of signature $n-1$ (see \cite{JV}). This would be enough
for timelike curves describing the trajectories of massive
particles, but if in addition we want to describe the trajectories
of massless particles,  we will consider a positive
two-homogeneous function $L$, rather than a one-homogeneous
function, ---since  for a Lorentzian metric $g$, $\sqrt{-g(v,v)}$
is not smooth when $v$ is lightlike\footnote{ There are some
non-trivial issues regarding smoothability. On the one hand,
mathematically, it is well known that all Finsler metrics smoothly
extendible to the 0-section must come from Riemannian metrics
\cite[Prop. 4.1]{Wa}. On the other, physically, there are
situations where it is natural to consider non-differentiable
directions, so that one may  allow this possibility explicitly
\cite{LPH}. However, we will not go through these questions in the
remainder (essentially, ``smooth'' might include some residual
non-differentiable points with no big harm).}. Summing up, we say
that a {\it (conic) Finsler spacetime} is a manifold endowed with
a conic pseudo-Finsler metric $L:A\rightarrow[0,+\infty)$ which
satisfies the following properties:
\begin{itemize}
\item[(i)] each $A_p:=A\cap T_pM$ is convex in $T_pM$ (i.e., the
segment in $T_pM$ connecting each two vectors $v_p, w_p\in A_p$ is
entirely contained in $A_p$), in particular, each $A_p$ must be
strictly included in a half-plane of $T_pM$, \item[(ii)] $A$ has a
smooth boundary in $TM\setminus \bf 0$, \item[(iii)] the extension
of $L$ as 0 to the closure $\bar A$ of $A$ is smooth at
$\hat{A}:=\bar A\setminus \bf 0$ (notice that the extension is
always continuous at $0$ by homogeneity, and this extension cannot
be smooth there even for a classical Finsler metric, except if it
comes from a Riemannian metric), and \item[(iv)] both the
fundamental tensor $g$ in \eqref{fundten} and its extension to the
points in $\partial A\setminus \bf 0$ have signature $n-1$.
\end{itemize}
 For practical purposes, we will consider always that $L$ is extended to $\bar A$ and $g$ to $T_v(T_pM)$ for all $p\in M$ and $v\in \bar A_p\setminus \bf 0$.
 In order to  make computations related to lightlike geodesics, one can also assume that $L$ is
  smoothly extended (in a non-unique way) on a neighborhood of $\hat{A}$. We will call $F=\sqrt{L}$
 the {\em Lorentz-Finsler metric} (defined on $\bar A$) of the Finsler spacetime. Observe that our notion of  Finsler spacetime is based on a
``conic'' element which is essentially present in  most of
previous literature.
 In particular:
\begin{enumerate}
\item  In \cite{Beem}, the function $L$ is defined in all $TM$ and
the fundamental tensor is assumed to have signature $n-1$. The
restriction to one connected component of the subset
$L^{-1}[0,+\infty)\setminus 0$ in every point $p\in M$ gives a
(conic) Finsler spacetime as defined here.   In \cite{Perlick06,
LPH}, the definition is essentially as in \cite{Beem} with the
opposite sign of $L$  (and an increasing attention to relax
differentiability).  The definition  in \cite{alemanes} is
somewhat  more involved: they introduce an $r$-homogeneous
function $L$, with $r\geq 2$ endowed with a fundamental tensor of
signature one and, then, they associate a homogeneous Finsler
function $F=\sqrt[r]{L}$. This allows   one  to deal with
non-differentiability in cases that extend  the lightlike vectors
 in our definition above.  In
\cite{Asanov}, the author   gives a definition analog to ours, but
he excludes the conditions on the boundary (so that one is not
worried about lightlike vectors).


 \item  Lorentz-Finsler metrics appear in the context of Lorentz violation.
 Here, one starts with a background Lorentz-Minkowski space, but
 the movement of particles or waves is governed by a Lagrangian with
 a Lorentz-Finsler behavior (see \cite{Russell, KoRu,KosRuTso} and
 references therein).
 Under natural hypotheses, the following  positively homogeneous
 function is found in \cite{kostelecky}:
\[F(v)= m \sqrt{-g_0(v,v)}+g_0(v,a)\pm\sqrt{g_0(v,b)^2-g_0(b,b)g_0(v,v)}, \]
where $g_0$ is the standard metric of Minkowski spacetime
$\R_1^4$, $v$ is any timelike vector and $a,b$ are two prescribed
vector fields in $\R_1^4$. The interplay between $F$ and the
background metric $g_0$ becomes important for physical
applications  (for example,  the group velocity of propagating
waves may exceed the light speed in $\R_1^4$). Even though $F$ can
be regarded as a Lorentz-Finsler metric, in some cases the
fundamental tensor of $F$ may be degenerate on some timelike
directions (when $b=0$, it is easy to compute the fundamental
tensor of $F$ \cite[Proposition 4.17]{JS}, and this becomes
degenerate when $F(v)=\sqrt{-g_0(v,v)}+g_0(v,a)=0$).

\item General physical theories of modified gravity (see
\cite{Vaca12} or the very recent articles \cite{SVaca, Vaca13} and
references therein) yield naturally Lorentz-Finsler spacetimes as
defined above. The exact domain $A$ of $L$ for general expressions
in Finsler-Cartan gravity (say, as in \cite[formula (2)]{Vaca13})
may be difficult to compute. However, one can compute   $A$ easily
in {\em General Very Special Relativity}, which possesses a simple
Finslerian line element introduced by Bogoslovsky (say, the
expression $F=(\sqrt{-g_0})^{(1-b)}\omega^b$ in \cite[formula
(I.1)]{KSS} becomes a Lorentz-Finsler metric with domain $A$ equal
to the intersection at each point of  the future time-cone of
$g_0$ and the half-space $\omega>0$) or in Randers-type
cosmologies as \cite[formula (3.1)]{BasStavr}.

\end{enumerate}


Notice that our definition of Finsler spacetime is  a
generalization  of the structure obtained in a classical spacetime
when one considers only its {\em future causal cones}. In the case
of Lorentzian metrics, however, the value of the metric on the
causal vectors is enough to determine the metric on all the
vectors and, moreover, the lightlike vectors are enough to
determine the metric up to a conformal factor. So, the cone
structure for Lorentzian metrics is equivalent to the conformal
structure. Nevertheless, this does not hold by any means in the
case of (conic) Finsler spacetimes: clearly two Lorentz Finsler
metrics $F, F'$ on $M$  with the same domain $\bar A$ may not be
equal up to a multiplicative function. However, the domain $A$ is
a {\em cone structure} and, thus, one can reconstruct all the
Causality Theory for Finsler spacetimes. Let us review this
briefly.

Given a Finsler spacetime $(M,L)$ a tangent  vector $v\in TM$ is
called (future-directed) timelike if $v\in A$, causal if $v\in
\hat{A}$, lightlike if $v\in \hat{A}\setminus A$ and null if
either $v$ is lightlike or the zero vector\footnote{Notice
that we have
 defined the Finsler metric only for ``future-directed'' causal
directions, i.e., those in  $\hat A$, and the past directions
could be consistently regarded as those in  $-\hat A$. It would be
natural then to extend $L$ to $-\hat A$ by requiring full
homogeneity ($L(v)=L(-v)$ for all $v\in -\hat V$), but this is not
necessary for the issues regarding Causality. However such an
extension (as well as possible further extension of $L$ to
non-causal vectors) would be important for the particular physical
theory one may be modelling.}. A (piecewise smooth) curve
$\gamma$ on $M$ is also called timelike, causal etc. depending on
the character of its velocity at all the points. If $p,q\in M$, we
say that  $p$ lies in the chronological (resp. causal) past of $q$
if there exists a (future-directed) timelike (resp. causal or
null) curve starting at $p$ and ending at $q$; in this case we
write $p\ll q$ (resp. $p\leq q$) and we also say that $q$ lies in
the chronological (resp. causal) future of $p$. The chronological
and causal futures of $p$, as well as its corresponding pasts, are
then defined formally as in the case of Lorentzian metrics:
$$\begin{array}{lr}

I^+(p)=\{q\in M:p\ll q\}, & J^+(p)=\{q\in M:p\leq q\}; \\
I^-(p)=\{q\in M:q\ll q\}, & J^-(p)=\{q\in M:q\leq p\}.

\end{array}$$
One says also that a second Lorentz Finsler metric $F'$ has cones
wider than $F$, denoted $F\prec F'$ when their corresponding
domains $A, A'$ satisfy $\hat{A} \subset A'$. With these
definitions, one can extend directly  the  {\em causal ladder} of
classical spacetimes (see for example \cite{MinSan}) to Finsler
spacetimes. We recall some of the steps of this ladder. A Finsler
spacetime will be called {\it chronological} (resp. {\it causal})
when it does not admit closed timelike (resp. causal) curves, {\it
stably causal} when there exists a causal  Lorentz-Lorentz metric
$F'$ with wider cones $F\prec F'$, and {\it globally hyperbolic}
when it is stably causal and the following property holds:
$J^+(p)\cap J^-(q)$ is compact for any $p,q\in M$. Trivially:
\[\text{chronological $\Leftarrow$ causal $\Leftarrow$ stably causal $\Leftarrow$ globally hyperbolic.}\]

Taking into account the case of spacetimes,
one realizes that there are more intermediate
levels of the ladder  as well as  many subtle properties and
relations among them to be considered for Finsler spacetimes.
However, we focus here just on a pair of them, with deep
implications for the global structure. To this aim, we  define for
a Finsler spacetime $(M,L)$:
\begin{itemize}
\item[(a)] a {\em spacelike hypersurface} is a smooth hypersurface
$S$  such that  no causal vector $v$ is tangent to $S$, \item[(b)]
a {\em (spacelike) Cauchy hypersurface} is a spacelike
hypersurface $S$ such that  any  causal curve which is
inextendible in a continuous way, intersects $S$ exactly once,
\item[(c)] a {\em temporal function} is a smooth function $t$ such
that $dt(v)>0$ for all (future-directed) causal vector $v$ (thus,
its levels $t=$constant are spacelike hypersurfaces) and,
\item[(d)] a {\em Cauchy temporal function} is a temporal function
such that all its levels  are Cauchy hypersurfaces.
\end{itemize}
In the case of classical spacetimes, results by Geroch \cite{Ge}
and Hawking \cite{Ha} obtained at a topological level,  plus their
improvements to the smooth and metric cases by Bernal and
S\'anchez \cite{BS1, BS2}, prove the equivalence between being
stably causal (resp, globally hyperbolic) and admitting a temporal
(resp. Cauchy temporal) function. By using different arguments
coming from KAM theory, these results were re-proved and extended
to general cone structures by Fathi and Siconolfi \cite{FS}. So,
as a consequence of the latter  we get the following result.

\begin{theorem}\label{tFS} Let $(M, F)$ be a Finsler spacetime.

(1) If $(M,F)$ is stably causal, then it admits a temporal
function (and, thus, it can be  globally foliated  by spacelike
hypersurfaces).

(2)  If $(M,F)$ is globally hyperbolic, then it admits a Cauchy hypersurface and  a Cauchy temporal function.

\end{theorem}
Recall, however, that the Finsler spacetime has elements that are not characterized by the cone structure. For example, one can define a {\em Finler separation} $d(p,q)\in [0,\infty]$ between any $p,q\in M$ by taking the supremum of the lengths of the (future-directed) causal curves from $p$ to $q$ ---in a close analogous to the {\em Lorentzian separation} or {\em distance} for spacetimes with Lorentzian metric.  The analogies and differences of  $d$ with the Lorentzian case (in particular, its interplay with Causality) deserve to be studied further.

 \begin{remark}
 One expects  that the converse of parts (1) and (2) of Theorem
\ref{tFS}, as well as most of classical causality theory, hold. To
achieve this,  typical Lorentzian techniques, as those regarding
limit curves or quasi-limits, must be adapted (recall the
Lorentzian proofs of the converse in \cite{Ge}, \cite[Sect.
6.6]{HE} or \cite[Th. 14.38, Cor. 14.39]{oneill}).  Moreover,
as the results in \cite{FS} work for very general cone structures,
Theorem \ref{tFS} can be applied to a more general class of
spacetimes that includes  those in \cite[Example 3.3]{Perlick06}, which
can be non-smooth in some timelike directions, and also the cone
structure provided by the example (2) above (see
\cite{kostelecky}).
 \end{remark}

\section{Stationary to Randers correspondence and beyond}\label{s4}
 From a
classical viewpoint,  a  correspondence between some purely
geometric elements of Lorentzian and Finslerian manifolds has been
developed recently. Such a correspondence  provides  a precise
description of certain objects in (classical) spacetimes in terms
of an associated Finsler space. The precise correspondence is
developed between a particular class of spacetimes, the {\em
standard complete-conformastationary} ones, or just {\em
stationary}, for short, and a precise class of Finsler manifolds,
the Randers ones (see \eqref{randers}). But some consequences and
techniques can be extrapolated to general Lorentzian and
Finslerian manifolds.

We start with a trivial observation. Probably the simplest examples of Lorentzian manifolds are the products $(\R\times M,g_L=-dt^2+\pi^\star g_0)$, where $\pi:\R\times M\rightarrow M$ and $t:\R\times M\rightarrow \R$ are the natural projections. The geometric properties of $(\R\times M, g_L)$ depend on those of
 $(M, g_0)$. In particular, each curve $(\R\supset) I\ni t\mapsto c(t)\in M$ parametrized with unit speed yields naturally two lightlike curves $t\rightarrow (\pm t,c(t))$ (future-directed with the sign $``+"$ and past-directed with $``- "$, for the natural time orientation of the spacetime), which are geodesics iff $c$ is a geodesic in $(M,g_0)$. A less trivial spacetime is obtained if we admit cross terms between the time and space parts (independent of the time $t$). This can be described by using a 1-form $\omega$ on $M$, namely, we consider the spacetime:
 \begin{equation}\label{estat}
 V=(\R\times M, g_L), \quad \quad g_L=-dt^2 + \pi^\star\omega \otimes dt + dt\otimes\pi^\star\omega +\pi^\star g_0.
 \end{equation}
Now, we introduce the following Finslerian  {\em Fermat metrics} associated to \eqref{estat}: $$F^\pm=\sqrt{g_0+\omega^2}\pm \omega,$$
notice that these are metrics of Randers type, and $F^-$ is the reverse Finsler metric of $F^+$, so, we will write simply $F$ for the latter and $\tilde F$ for the former. If we consider curves $c^+$ and $c^-$ in $M$ which are unit for $F$ and $\tilde F$ resp.,  the curves in the spacetime $t\rightarrow (\pm t,c^\pm (t))$ are again (future or past directed) lightlike curves in  $V$, and each one is a geodesic  up to parametrization  iff so is the corresponding original curve $c^+$ or $c^-$ (see the details in \cite{CJM11} or \cite{CJS}). This suggests the possibility  of describing the properties related to the Causality and conformal structure in \eqref{estat} in terms of the geometry of the corresponding Randers space $(M,F)$.

It is worth pointing out that, as only conformally invariant
properties will be taken into account, the class of spacetimes to
be considered includes those conformal to  \eqref{estat}. This
class of spacetimes can be characterized intrinsically as those
which are distinguishing (a causality condition less restrictive
than strong causality) and admit a complete timelike conformal
vector field. In fact, the conformal change of the metric $g$ by
$-g/g(K,K)$, plus the result in \cite{JS}, allows  one to find the
expression \eqref{estat}, which  defines a {\em normalized
standard stationary spacetime}. So, one will find a correspondence
between the conformal properties of the elements in this class of
metrics $(V,g_L)$ and the geometric properties of Randers spaces
$(M,F)$. This has been carried out at different levels (see
\cite{BJ11,CGS,CJM10,CJM11,CJS,FHSmem,GHWW09,JLP} or \cite{J} for
a review), and we will focus here in three of them, with clear
physical applications.

\subsection{Causal structure}

Nicely,  Fermat metrics allow  one to determine the chronological
and causal future and past of any point in the stationary
spacetime \eqref{estat}. For example, if we take $(0,p)\in
\{0\}\times M$ then the intersection of $I^+(0,p)$ with the slice
$t=t_0$ of the spacetime is equal to $\{t_0\}\times B^+(p,t_0)$,
where $B^+(p,t_0)$ is the open ball of center $p$ and radius $t_0$
for the Fermat metric $F(=F^+)$.  From these considerations, one
can  describe in a precise way the causal structure of the
spacetime. Concretely, one has \cite{CJS}:

\begin{theorem}\label{t1}
A stationary spacetime \eqref{estat} is always causally continuous and it is
\begin{itemize}
\item Causally simple (i.e., it is causal, and the causal futures
and pasts $J^\pm(p)$ are closed) if and only if $(M,F)$ is convex,
i.e. any pair $(p,q)$ in $M$ can be connected by means of an
$F^+$-geodesic of length equal to the Finslerian distance
$d_F(p,q)$. \item Globally hyperbolic if and only if the closed
balls for the symmetrized distance of $d_F$ are compact. \item
Globally hyperbolic with slices $t=$constant that are Cauchy
hypersurfaces if and only if $d_F$ is forward and backward
complete.
\end{itemize}

\end{theorem}
Other causal elements which are described naturally with the Finslerian elements are the Cauchy horizons and developments. For example, given a subset $A\subset \{0\}\times M$, its future Cauchy horizon $H^+(A)$ is the graph, in $\R\times M$ of the function which maps each $y\in \bar A$ to $d_F(M\setminus A, y)$, i.e., the $d_F$-distance from the complement of $A$ to $y$.

\begin{remark}
The results in  stationary spacetimes  can also be translated to results in Randers metrics, sometimes generalizable to any Finsler manifold. Among them, we point out (see \cite{CJS}):

\begin{itemize}
\item The theorem above suggests that the compactness of the closed symmetrized balls of $d_F$, which is a weaker hypothesis than the commonly used (forward or backward) completeness of $d_F$, can substitute the last hypothesis in  many results. And, in fact, this is the case in  classical theorems of Finsler Geometry, such as  Myers' theorem or the sphere theorem. 

\item From the known fact that any globally hyperbolic spacetime admits a smooth spacelike Cauchy hypersurface, one can deduce that any Randers metric $R$ with compact symmetrized balls admits a trivial projective change ($R\rightarrow R+df$) such that the corresponding new Randers  metric has {\em the same pregeodesics as $R$, and it is forward and backward complete}. This result has been extended by Matveev \cite{Ma} for any Finsler metric.

\item  The results on horizons (a substantial topic in Lorentzian Geometry) in stationary  spacetimes, yield directly results on Randers spaces. For example, the translation of a well-known result by Beem and Krolak \cite{BK} yields the following property for any subset $A$ of a Randers manifold $(M,R)$: {\em $p\in M$ is a differentiable point of the distance from $p$ if and only if it can be crossed by exactly one minimizing segment.} This result has been extended to any Finslerian manifold by Sabau and Tanaka \cite{ST}.
\end{itemize}
\end{remark}

\subsection{Visibility and gravitational lensing}

Assume that in our spacetime \eqref{estat}, a point $w$ represents
an event and a line $l$, obtained as an integral curve of
$\partial_t$, represents the trajectory of a stellar object. We
wonder if there    exist lightlike geodesics from $l$ to $w$
(i.e., whether an observer at $w$ can see the object $l$) and, in
this case, if there are many  of such geodesics (i.e., the
existence of a lens effect such that $l$ is seen in two different
directions). This situation is applicable to cosmological models
such as Friedmann-Lemaitre-Robertson-Walker ones, as they are
conformally stationary ($l$ would represent a ``comoving
observer'' of the model). The problem becomes more realistic if we
{ choose} some (open) region $\R\times D\subset \R\times M$ which
contains $w$ and $l$, and  search for geodesics contained in this
region.

These problems are related to the convexity of $D$ for the Fermat metric (in the sense of the last subsection) which turns to be related to the convexity of its boundary $\partial D$. There are several different notions for the latter, as the {\em local} and {\em infinitesimal convexity} and, as shown in \cite{BCGS}, they are equivalent to its {\em geometric convexity}. The latter means that, given any pair of points  of $D$, any geodesic connecting them and contained in the closure $\bar D$, must be entirely  contained in $D$. The (geometric) convexity of $D$ turns out equivalent to the {\em light (geometric) convexity} of the boundary of $\R\times D (\subset \R\times M)$ (i.e., the property of  convexity holds when restricted to lightlike geodesics) and, finally, this is equivalent to the question of  existence of connecting geodesics, yielding \cite{CGS}:

\begin{theorem}
Assume that the closed balls of $\bar D$ computed for the restriction  of the symmetrized distance $d_F^s$
are compact (which happens for example, if the intersections with $\bar D$ of the closed symmetrized balls in $M$ are compact or, simply, if $d_F$ is complete on all $M$). Then, the following assertions are equivalent:
\begin{romanlist}[(iii)]
\item $(D\times\R,g_L)$ is causally simple (i.e. $(D,F)$ is convex, Theorem \ref{t1}).

 \item $\partial D$ is convex for the Fermat metric $F$.

 \item $ \R \times\partial D$ is light-convex for the Lorentzian metric $g_L$.

\item Any point $w  = (t_p,p) \in \R \times D$ and any line
$l_q:=\{ ( \tau,q ) \in \R \times D: \tau \in \R\}$, with $p\neq
q$, can be joined in $\R\times D$ by means of a future-directed
lightlike geodesic $z (s) = (  t (s),x (s)), s\in [a,b]$, which minimizes the
(future) arrival time $T=t(b)-t(a)$  (i.e., such that
$x$ minimizes the $F$-distance in $D$ from  $p$ to $q$).

\item Idem to the previous property but replacing ``future'' by
``past'' and ``$F$-distance'' by (the reverse) ``$\tilde{
F}$-distance''.

\end{romanlist}
\end{theorem}

\begin{remark} The previous result solves the question of existence of connecting lightlike geodesics. The question of multiplicity has some possibilities. First, the existence of a conjugate point of the lightlike geodesic $z$ at $w$. This is equivalent to the existence of a conjugate point at $p$ for its projection $x$ \cite[Theorem 13]{CJM10}, and it is regarded as trivial.  Second, the non-triviality of the topology of $D$ may yield a {\em topological lensing}. In fact, one can prove that, whenever $D$ is not contractible, infinitely many connecting lightlike geodesics (with diverging arrival times) will exist.

The results can be also extended to the case of timelike geodesics, prescribing its length (i.e., the lifetime of the massive particles represented by such geodesics). The idea relies on a  reduction of the problem to the lightlike one, by considering an extra spacelike dimension, see \cite[Section 4.3]{CJS} and \cite[Section 5.2]{CGS}.
\end{remark}

\subsection{Causal boundaries and further questions}

The studied  {\em stationary to Randers correspondence} can be also applied to the study of the causal boundary of the spacetime in terms of Finsler elements. We recall that, in Mathematical Relativity, the Penrose conformal boundary is commonly used, in spite of the fact that it is not an intrinsic construction, and  there are problems to ensure its existence or uniqueness.  The causal boundary is a conformally invariant alternative,  which is intrinsic and can be constructed systematically in any strongly causal spacetime (see \cite{FHSatmp} for a comprehensive study of this boundary). The computation of the causal boundary and completion of  a stationary spacetime \eqref{estat} has been carried out in full generality in \cite{FHSmem}. It must be emphasized that this boundary has motivated the definition of a new {\em Busemann boundary} in any Finslerian manifold. Even more, this has stimulated the study of further properties of the Gromov boundary in both, Riemannian and Finlerian manifolds.

It is also worth mentioning two further topics of current interest. The first one is the correspondence
 at the level of curvatures between the Weyl tensor of the stationary spacetime and the flag curvature of the Randers space, see \cite{GHWW09}. The second one is the possibility to extend the studied correspondence to the case of spacetimes with a nonvasnishing complete Killing vector field, whose causal type may change from timelike to lightlike and spacelike. Such metrics are related to a generalization of the Zermelo metrics which includes the possibility of a strong wind (i.e., with a speed higher than the one which can be reached by the engine of the ship), see \cite{CJSpro}.
 \section{Conclusions}
As a first goal, we have introduced a very general notion of
Finsler spacetime which includes many previous ones. This notion
emphasizes the role of the cone structure, showing that many
properties of (classical, relativistic)  spacetimes are
potentially transplantable to this setting ---in particular, this
holds for a major Lorentzian result such as Theorem 1. Remarkably,
in relativistic spacetimes the cone structure is equivalent to the
conformal structure and, thus, it is governed by the Weyl tensor.
However, the metrics $F_1, F_2$ of two Finsler spacetimes with the
same cone structure are not by any means necessarily conformal,
and one is lead to deal directly with the cone structure. This
cone structure would remain then at a more basic level than any
possible Weyl-type tensor  and, thus, also at a much basic level
than the Cartan or Chern connections (of course, such connections
will be relevant in order to describe concrete physical effects).

The importance on the choice of the connection is stressed because
a Finsler metric admits several associated connections.  Among
them, one has the classical  Cartan connection, the Berwald
connection, the Chern-Rund connection and the Hashiguchi
connection (see \cite{SzLo} or \cite[\S  9.3]{SzLoKe24}) .  Thus, for example,  two Finsler
metrics with the same Chern connection may have different Cartan
connections.
 As a consequence,  the way to
describe physical phenomena by using these connections may be 
non-equivalent and  even controversial \cite{ChLi,Vaca10}.
Nevertheless, as emphasized here, all the possible physical causal
influences will depend exclusively on the cone structure, at least
at a non-quantum level.

As a second goal we have shown a different type of equivalence
between the geometry of any Randers space and the conformal
geometry of classical standard stationary spacetimes. This
equivalence leads to obtain very precise results which relate: (a)
the properties of distances on the Finslerian side with the
Causality on the side of the spacetime, (b) the properties related
to Finslerian {\em convexity} wth the properties related to
relativistic observability and lensing, (c) accurate relations
between ideal Finslerian  boundary  and the causal boundary (and,
implicitly then, the conformal one) of spacetimes, or (d)
interesting relations between Finslerian curvature and Weyl
curvature \cite{GHWW09}. This equivalence has suggested new
results in both, Finslerian and Lorentzian geometries. Moreover, a
forthcoming work \cite{CJSpro}  will show its extendability and
applicability to new fields ---namely, an extension of classical
Finsler metrics which allows one to model strong winds in Zermelo
navigation problem is introduced, and its geometrical properties
are shown equivalent to the conformal ones of a class of
spacetimes wider than standard stationary ones, which allow the
possibility of horizons and black holes.

\section*{Acknowledgments}
 Comments by  Kostelecki, Laemmerzahl, Perlick,  Szilasi  and the referee
are warmly acknowledged. Both authors are  partially supported by
the Grant P09-FQM-4496 (J. Andaluc\'{\i}a) with FEDER funds. The
first-named author is also partially supported by
MINECO-FEDER
project MTM2012-34037 and
Fundaci\'{o}n S\'{e}neca project 04540/GERM/06, Spain. This
research is a result of the activity developed within the
framework of the Programme in Support of Excellence Groups of the
Regi\'{o}n de Murcia, Spain, by Fundaci\'{o}n S\'{e}neca, Regional
Agency for Science and Technology (Regional Plan for Science and
Technology 2007-2010).  The second-named author acknowledges
the support of IHES, Bures sur-Yvette, France, for a three months
stay in 2013 (which included the period of celebration of the
meeting IFWGP); he is also partially supported by MTM2010--18099
(MICINN-FEDER).

\end{document}